\documentclass[11pt]{amsart}
\usepackage{amsfonts}
\usepackage{color}
\usepackage{graphicx}

\theoremstyle{plain}
\begingroup
\newtheorem{theorem}{Theorem}[section]

\endgroup

\theoremstyle{definition}
\begingroup
\newtheorem{defn}[theorem]{Definition}
\newtheorem{remark}[theorem]{Remark}
\newtheorem{rmk}[theorem]{Remark}

\endgroup

\theoremstyle{remark}
\begingroup
\endgroup

\mathchardef\emptyset="001F

\numberwithin{equation}{section}

\newcommand{\e}{\varepsilon}
\newcommand{\Om}{\Omega}
\newcommand\io{\int_\Om}
\newcommand\xo{x_0}

\newcommand{\R}{{\mathbb R}}

\newcommand{\rn}{{\R}^N}
\newcommand{\rt}{{\R}^2}

\newcommand{\hn}{{\mathcal H}^{N-1}}
\newcommand{\huno}{{\mathcal H}^{1}}

\newcommand\ee{\end{equation}}
\newcommand\be{\begin{equation}}
\newcommand\wc{\rightharpoonup}
\newcommand\ep{\varepsilon}
\newcommand\uh{\hat{u}}
\newcommand\kh{\hat{K}}

\newcommand\nb{\partial_N\Omega}
\newcommand\db{\partial_D\Omega}

\title[A local variational principle for fracture]
{A local variational principle for fracture}

\author[Christopher J. Larsen]
{Christopher J. Larsen}
\address[Christopher J. Larsen]{Department of Mathematical Sciences,
Worcester Polytechnic Institute,
100 Institute Road,
Worcester MA 01609-2280,
USA}
\email[]{cjlarsen@wpi.edu}

\begin{document}

\begin{abstract}
The seminal paper of Francfort and Marigo \cite{fm} introduced a variational formulation for Griffith fracture \cite{g} that has resulted in substantial theoretical and practical progress in modeling and simulating fracture.  In particular, it led to the phase-field approximation proposed in \cite{bfm1}, which has been widely implemented.  However, the formulation in \cite{fm} is known to have limitations, including its inability to treat applied loads and its reliance on global minimization. 
In addition, the phase-field model \cite{bfm1} and its extensions, as implemented, are not generally approximations of the global minimizers in \cite{fm}.  In this paper, we show that there is a local variational principle satisfied by global and local minimizers of the energy introduced in \cite{fm}, which is compatible with loads, and which is a generalization of the stress intensity factor.  We use this principle to reformulate variational fracture, including formulations that, for the first time, can include all forms of applied loads.  We conclude by showing the connection between phase-field models, as implemented, and our formulations.
\end{abstract}
\maketitle
{\small
 
\bigskip
\keywords{\noindent {\bf Keywords:} {Variational fracture, phase-field fracture, $\Gamma$-convergence.}

\bigskip
\bigskip

\section{Introduction}

The variational formulation of Griffith fracture introduced by Francfort and Marigo \cite{fm} has led to major advances in the mathematical study and simulation of fracture.  A central feature of that formulation is its mathematical robustness -- it can be stated with no assumptions at all on crack geometry or regularity, such as the existence of crack tips.  From a mathematical point of view, this is very advantageous, allowing existence proofs (which are missing for models based on stress intensity factors) as well as rigorous approximation results, such as those based on phase-fields.

Benefits of this formulation are evident when implementing phase-field fracture, since ``cracks" grow where elastic energy is  sufficiently large, again with no assumption whatsoever about the geometry of the ``cracks".  On the other hand, this locality, while physically desirable, is not present in the original formulation \cite{fm}, since the formulation there is based on global minimization (and implementations of phase-field fracture do not involve global minimization).  A natural and important question is then, what are these phase-field formulations approximating?  The limiting sharp interface model must not be based on global minimality, but instead some kind of local property.

An additional limitation of global minimization is that it is incompatible with applied body loads.  Our goal is to find  a  local variational principle that works seamlessly with loads, while having the same mathematical robustness (no regularity or geometric assumptions on the crack) as \cite{fm}.  

Another reason to seek a local formulation  is to allow nucleation based on strength \cite{kfl}.  There, growth of existing (large) cracks is governed by Griffith's criterion, whereas nucleation of new cracks is based on completely separate ``strength" criteria.  Global minimality as in \cite{fm} can force nucleation, independent of the specified nucleation criteria.  Therefore a local Griffith criterion, governing only growth of existing cracks, is desirable.

Below, we describe a local variational principle --  a necessary condition for Griffith stability -- which naturally works with all types of applied loads.  Furthermore, it is weak in the sense we want, in that there is no requirement of any regularity of cracks, such as the existence of crack tips.  
  In essence it allows us to specify exactly what energy competes with fracture surface energy, separately from specifying the energy that determines the elastic equilibrium.  Furthermore, it seems to be the first necessary condition for Griffith stability that requires no regularity of the crack.

This paper is organized as follows.  In Section \ref{limits} below, we describe in more detail limitations of current sharp interface variational models. Interestingly, these problems are largely missing from phase-field models, as implemented, and we describe these features of phase-field models in Section \ref{fpf}.  In Section \ref{fpvp} we state the proposed variational principle, {\em local Griffith stability}.  We then give precise variational formulations for different models, both static and quasi-static, with and without loads, and even define static crack stability for dynamic fracture.  We also describe how this new stability definition might be used to formulate models when only part of the elastic energy competes with Griffith surface energy, as proposed for viscoelastic materials in \cite{sl}.

Section \ref{stress} discusses the relationship between stability based on stress intensity factors, and the stability we propose here.  We show their consistency, and that our approach is more general, in that no extra assumptions, such as regularity or crack geometry, are required.  Then in
Section \ref{mgs} we describe why  energy minimizers satisfy  local Griffith stability.  
Finally, in Section \ref{gcpf}, we discuss the connection between our formulations and phase-field models.  In particular, we describe the role $\Gamma$-convergence plays in relating rescaled phase-field models to our variational principle.

\subsection{Limitations of existing sharp interface models: the problem with loads}\label{limits}

Before we explain the problem with loads, we first need to be more precise about the formulation in \cite{fm} without loads. Considering for simplicity scalar displacements $u$ on $\Om\subset\rn$  with  the simplest elastic energy density $\frac12 |\nabla u|^2$, Griffith constant $G_c$,  and crack sets $K$, the total energy of a displacement-crack pair $(u,K)$, with $u$ allowed to be discontinuous across $K$, is given by
\[
E^G(u,K):=\frac12 \io |\nabla u|^2 dx + G_c \hn(K),
\]
where $\hn(K)$ is the $N-1$ dimensional Haudorff measure (surface area) of the crack $K\subset\bar{\Om}$, and $u\in H^1(\Om\setminus K)$ (a slight abuse of notation explained below). 
Minimizing this energy (globally or locally) subject to a given displacement boundary condition produces a displacement $u$ and crack set $K$ satisfying $\Delta u=0$ in $\Omega\setminus K$ with normal derivative $\partial_\nu u=0$ on $K$.  In addition, the crack $K$ results from direct competition between elastic and surface energy, producing a form of Griffith stability:
\[
E^G(u,K)\leq E^G(v,\kappa)
\]
for all pairs $(v,\kappa)$ satisfying $v=u$ on $\partial \Om$ and $\kappa\supset K$.  If $(u,K)$ are obtained by global minimization, then this inequality holds for all such $(v,\kappa)$, and if they are obtained by local minimization, then the inequality will hold if $v$ is sufficiently close to $u$.

We now turn to the problem of including   boundary and body loads in variational fracture.  With boundary load $g$ on part of $\partial\Om$ denoted $\nb$, and displacement boundary condition zero, for example, on the remainder of the boundary $\db$, it would seem that we should minimize
\[E_{Load}(u, K):= \frac12 \io |\nabla u|^2 dx-\int_{\partial_N \Omega}gu\; ds +G_c \hn(K)\]
over $u\in H^1(\Omega\setminus K)$ with $u=0$ on $\db$.
Minimizing this energy is easily seen to be impossible (except in the trivial case $g\equiv 0$).  The idea is, we can choose a part of $\nb$ in which the average of $g$ is not zero, and if we create a crack disconnecting that part of the boundary from $\db$, the second term in the energy can be sent $-\infty$ with a controlled cost in the rest of the energy, so the total energy goes to $-\infty$.   Essentially the same issue occurs with body loads.

There are additional issues with loads, besides nonexistence, which we describe in \cite{la} and do not repeat here.  
But these issues are not present if we want the crack surface energy to compete with elastic energy the same way it does with displacement boundary conditions, as we advocated in \cite{la}.  There, we showed that this principle can be implemented with boundary loads -- instead of trying to minimize a single energy, we can find $(u,K)$ simultaneously minimizing two different energies:
\[
E_K(v):=\frac12 \io |\nabla v|^2 dx -\int_{\nb} gv \;ds
\]
over displacements $v$ with the same crack and with $v=u$ on $\db$, and
\[
E^G(v,\kappa)=\frac12 \io |\nabla v|^2 dx + G_c \hn(\kappa)
\]
over displacement-crack pairs  with $v\in H^1(\Om\setminus \kappa)$, $v=u$ on all of $\partial \Om$, even $\nb$, and for $\kappa\supset K$.  We then get the elastic plus load equilibrium $u$ we want, but the crack results only from competition between elastic and surface energy (of course, the crack indirectly depends on the load, but only by the dependence of $u$ on the load, and not direct energy competition between the load and surface area).  This formulation, and a method for showing existence, are developed in \cite{la}.

Body loads are much worse. The total energy is now
\[
\frac12 \io |\nabla u|^2 dx -\io fu \;dx + G_c \hn(K),
\]
and there is a big problem with global and local minimizers, in addition to corresponding issues to those described above for boundary loads.  
Splitting this energy into two energies, as with boundary loads, cannot work, since $u$ cannot simultaneously minimize
\[
\frac12 \io |\nabla u|^2 dx -\io fu \;dx
\]
and
\[
\frac12 \io |\nabla u|^2 dx + G_c \hn(K),
\]
as minimizing the first implies $\Delta u=f$ in $\Omega\setminus K$, while minimizing the second gives $\Delta u=0$ in $\Omega\setminus K$.  So the question for body loads is, how can we produce $u$ minimizing the first energy, while having the crack in equilibrium only with elastic energy?  Below, we show how this can be done with a local variational principle.

A related question is raised by \cite{sl}, where they propose that for certain materials, the displacement $u$ should minimize an elastic energy of the form
\[
 \io W(\nabla u)  dx, 
\]
with $W=W^{Eq} + W^{NEq}$, but the elastic-fracture equilibrium should come from minimizing just
\[
 \io W^{Eq}(\nabla u) dx + G_c \hn(K).
\]
So here, too, the variational principle for determining $u$ (given $K$) cannot be combined with the variational principle for elastic-fracture equilibrium into one energy.  However, we will consider below splitting this into two variational problems, one for determining the equilibrium displacement, and the other a local principle that sees only $W^{Eq}$ and surface energy.

\subsection{Features of phase-field models, as implemented}\label{fpf}
In principle, none of the above are problems for phase-field fracture, {\em as implemented}.  While the elastic plus surface energy from \cite{bfm1} is
\[
E_\delta(u,v):=\frac12\io (\eta_\delta + v^2)|\nabla u|^2 dx + G_c\left( \frac1{4\delta} \io (1-v)^2 dx + \delta \io |\nabla v|^2 dx  \right),
\]
in practice, for static and quasi-static settings, this energy is effectively split into two energies, which are separately (alternately) minimized: with $v$ fixed, $u$ minimizes
\[
E_v(u):=\frac12\io (\eta_\delta + v^2)|\nabla u|^2 dx 
\]
and with $u$ fixed, $v$ minimizes
\[
E_u(v):=\frac12\io  v^2 |\nabla u|^2 dx + G_c\left( \frac1{4\delta} \io (1-v)^2 dx + \delta \io |\nabla v|^2 dx  \right).
\]
The first energy $E_v$ includes only the part of $E_\delta$ that depends on $u$, and the second, $E_u$, on $v$.  Notice that the latter naturally includes the elastic and surface energy, and minimizing in $v$ means that the crack advances when the reduction in the elastic energy compensates for the increase in surface energy, {\em independently of the rule governing $u$}.  For static and quasi-static fracture, $u$ is found by minimizing the first energy, $E_v$, given the updated $v$.  $E_u$ and $E_v$ are then alternately minimized until convergence.  For dynamics, elastodynamics is used to find $u$, where the stiffness is $\eta_\delta + v^2$ (\cite{blr}), and $v$ is found by minimizing $E_u$.

These procedures are easily modified to include loads.  For example, with body loads, it is enough to change the first energy only, and for statics and quasi-statics, minimize in $u$:
\[
E^L_v(u):=\frac12\io (\eta_\delta + v^2)|\nabla u|^2 dx -\io fu dx
\]
so that the minimal $u$ is in elastic equilibrium subject to the body load $f$, and there is no direct competition between the surface energy and the work done by the load.  Minimizing in $v$ is done just as before, with elastic energy alone competing with surface energy.  Boundary loads can be handled similarly (even in sharp interface fracture, see \cite{la}).  In general, whatever principle determines $u$, minimizing $E_u$ in $v$ results in the crack advancing only if the reduction in elastic energy alone is sufficient to justify the increase in surface energy.

This extends even to the situation in \cite{sl} described above, where for certain viscoelastic materials there are two parts to the elastic energy, and only one competes with surface energy: the equilibrium displacement is found by minimizing the total stored elastic energy
\[
E_v(u):=\io (\eta_\delta + v^2)W(\nabla u) dx 
\]
with $W=W^{Eq} + W^{NEq}$, but only $W^{Eq}$ competes with surface energy: $v$ minimizes
\[
E_u(v):=\io (\eta_\delta + v^2)W^{Eq}(\nabla u) dx + G_c\left( \frac1{4\delta} \io (1-v)^2 dx + \delta \io |\nabla v|^2 dx  \right).
\]
$u$ is an equilibrium displacement, with elastic energy $W$, but the crack advances based only on whether the reduction in the $W^{Eq}$ elastic energy compensates for increased surface energy.  Here, we can expect that alternate minimization produces $(u,v)$ with $u$ minimizing $E_v$ and $v$ minimizing $E_u$.

In addition, we emphasize that these phase-field formulations have the physically desirable property that they are spatially local: cracks advance ($v$ decreases) only if the elastic energy is large enough {\em where the crack advances}.

With the exception of boundary loads, none of the above can be done with existing sharp interface variational formulations.

\section{The local variational principle}\label{fpvp}

The main definitions we introduce are {\em {local Griffith stability}}, Definition \ref{lgs} below, and  {\em {elastic-Griffith stability}}, Definition \ref{maindefstatic}.  Local Griffith stability is a pointwise stability condition, in the same way that $\Delta u=f$ is a pointwise balance condition.  It is naturally applied at crack tips, and is a generalization of stress intensity factors, with the benefit that no crack tips, and no crack geometry at all, are necessary in order to state the stability condition.   

Before stating the condition, we first briefly discuss notation.  Mathematically, the appropriate settings for ``sharp interface" fracture are rather technical spaces based on $SBV$ or $SBD$.  For our purposes, we do not need to go into these technical details, and so we will use more intuitive notation which is only slightly an abuse.  We first recall that if a crack set $K$ is closed in $\Omega$, then an admissible displacement $u$, which can be discontinuous across $K$, is in the space $H^1(\Omega\setminus K)$. Below we will use this notation even when $K$ is not closed.

Second, if $u$ and $v$ are both in $H^1(\Omega\setminus K)$ with the same trace on $\partial \Omega$ (which is a technical form of boundary values), we say $v\in H_u^1(\Omega\setminus K)$, which is the same as saying $u-v\in H_0^1(\Omega\setminus K)$.  For simplicity,  we will usually take $\Omega\subset\rt$, even when the corresponding result holds in $\rn$.  When there is something special about $\rt$, we will note it.

 We say $u\in H_{loc}^1(\rt\setminus K)$ if for every bounded open set $U\subset\rt$, we have $u\in H^1(U\setminus K)$ (the point being that integrals over all of $\rt$ might not be finite).  Finally, when integrating over $\Omega\setminus K$, we usually write the integral over all of $\Omega$, since the $K$ we consider have Lebesgue measure zero.

We now turn to the proposed variational principle.
The idea is that under a certain rescaling, first introduced in \cite{bo}, blow-up limits of stable states must have a global minimality property.  We will save the details of taking limits, and why they should have this minimality, for Section \ref{mgs}, but we state the rescaling and minimality now.  Given a displacement-crack pair $(u,K)$ with $u\in H^1(\Om\setminus K)$, for $x_0\in \Om$ fixed, and for each $\e>0$, we set

\begin{figure}[h]
\includegraphics[scale=0.35]{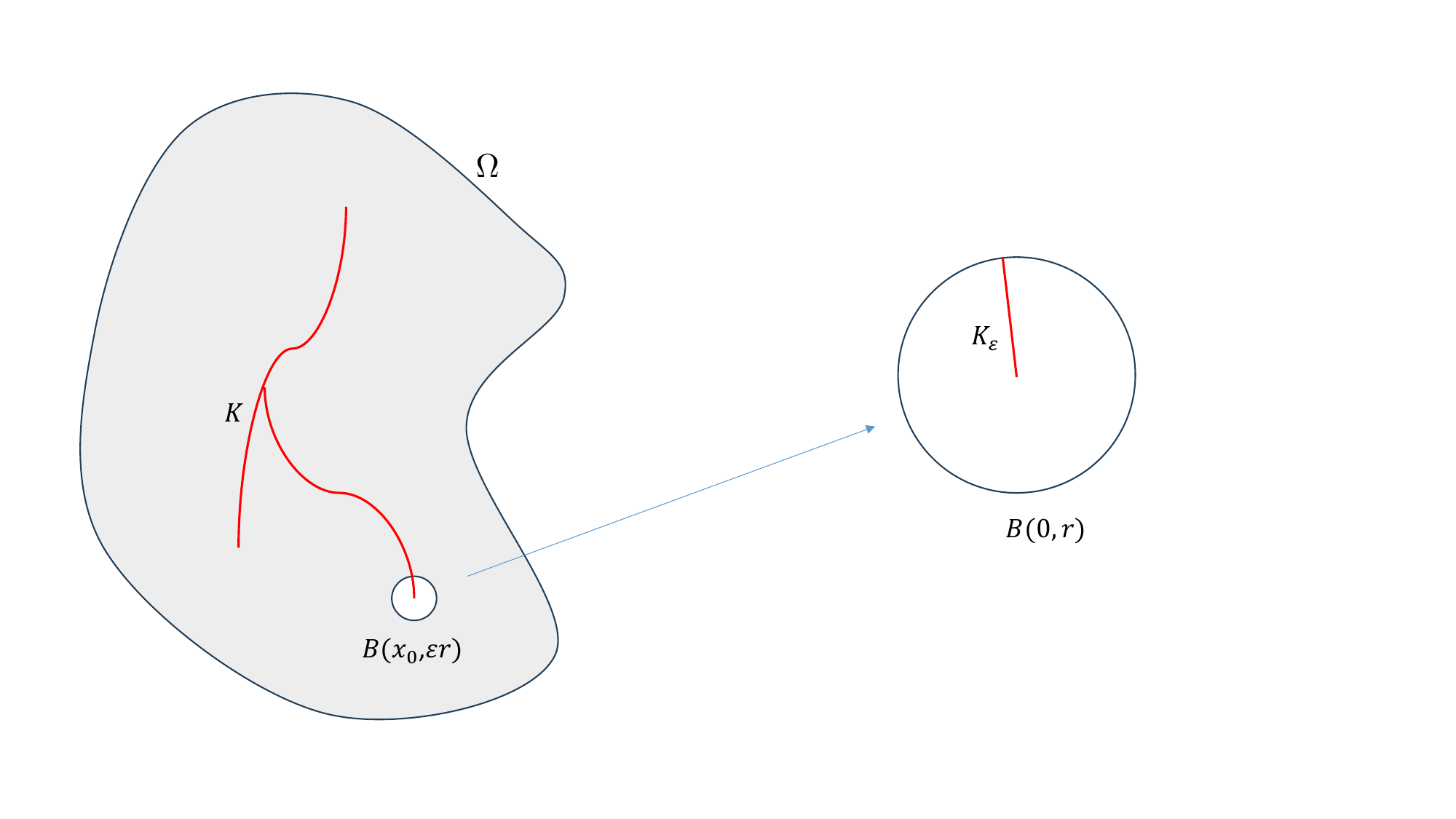}
\caption{Rescaling $K$ at $x_0$}
\label{blo}
\end{figure}

\begin{equation}\label{bu}
u_\e(x):=\ep^{-\frac12}[u(x_0+\ep x)-u(x_0)]\\
\mbox{ and }K_\e:=\e^{-1}(K-x_0),
\end{equation}
see Figure \ref{blo}.

If a pair $(\uh, \kh)$ is the limit as $\ep \rightarrow 0$ of the above rescaled pairs (or the limit of a subsequence), then due to the rescaling, both $\uh$ and $ \kh$ would be defined on all of $\rt$, and so $\uh\in H_{loc}^1(\rt\setminus \kh)$ and $\kh\subset \rt$.  The minimality that we will consider for this pair is

\begin{defn}[Global Griffith stability]\label{gs}
A pair $(\hat{u},\hat{K})$ with $\uh\in H_{loc}^1(\rt\setminus \kh)$ is {\em{globally Griffith-stable}} if for every $r>0$, it minimizes
\[E_r^G(w,\kappa) := \frac12 \int_{B(0,r)}|\nabla w|^2 dx +G_c\huno(\kappa \cap B(0,r))
\]
over pairs $(w,\kappa)$ satisfying $w\in H_{\uh}^1(B(0,r)\setminus \kappa)$ and $\kappa \supset \kh\cap B(0,r)$.  Here, $\huno$ is the one-dimensional Hausdorff measure.

\end{defn}

Our main definition just requires all such blow-up limits to have this minimality:

\begin{defn}[Local Griffith stability]\label{lgs}
$(u,K)$ with $u\in H^1(\Om\setminus K)$ is {\em locally Griffith stable} if for every $x_0\in\Om$, every blow-up limit $(\hat{u},\kh)$ of  $(u_{\e},K_{\e})$  is globally Griffith stable.

\end{defn}

 The convergence to blow-up limits is a bit technical, and we discuss this in more detail below, in Remark \ref{rmk}.

We now describe how this local stability can be used to define stable states and evolutions.  For the simplest setting, static fracture with no loads, we will  combine the usual elastic equilibrium condition with the new local Griffith stability condition.  Later we can make a very natural modification to include loads.  The formulation for quasi-static fracture with loads is a natural further modification.

\begin{defn}[elastic-Griffith stability]\label{maindefstatic}

We say $(u,K)$ with $u\in H^1(\Om\setminus K)$ is elastic-Griffith stable if:
\begin{enumerate}
\item $u$ minimizes
\[ w\mapsto \frac12\io |\nabla w|^2 dx \]
over $w\in H^1(\Omega\setminus K)$ with $w=u$ on $\partial \Omega$, i.e., 
\[\Delta u=0\mbox{ in }\Omega\setminus K,\;   \partial_\nu u=0\mbox{ on }K,\]
and

\item $(u,K)$ is locally Griffith stable (Definition \ref{lgs}).
\end{enumerate}

\end{defn}
Of course, this stability can be posed with specified boundary values, but we leave that out here since it is unnecessary for defining stable states.

Body and boundary loads are easily added:

\begin{defn}[elastic-Griffith stability with loads]\label{maindefstaticbody}

We say $(u,K)$ with $u\in H^1(\Om\setminus K)$ is elastic-Griffith stable with body load $f$ and boundary load $g$ (applied to part of the boundary $\nb$) if:
\begin{enumerate}
\item $u$ minimizes
\[ w\mapsto \frac12\io |\nabla w|^2 dx -\io fw \;dx -\int_{\nb} gw \;ds\]
over $w\in H^1(\Omega\setminus K)$ with $w=u$ on $\partial_D \Omega:=\partial \Om \setminus \nb$, i.e., 
\[\Delta u=f\mbox{ in }\Omega\setminus K, \;  \partial_\nu u=0\mbox{ on }K, \mbox{ and }\;\partial_\nu u=g\mbox{ on }\nb,\]
and

\item $(u,K)$ is locally Griffith stable.

\end{enumerate}

\end{defn}
Note as above that global and local minimizers of variational fracture with boundary loads (see \cite{la}) satisfy this definition, as do local minimizers with body loads (see Section \ref{mgs}).  Note also that if $f=0$ and $g=0$, Definition \ref{maindefstaticbody} reduces to Definition \ref{maindefstatic}, as we would want.

Next we have the general quasi-static formulation:

\begin{defn}[Quasi-static elastic-Griffith stable evolution with loads]\label{maindef}

We say $(u,K)$ is a quasi-static evolution on $[0,T]$ with body load $f$ and boundary load $g$ if it satisfies the following:
\begin{enumerate}
\item for every $t\in [0,T]$, $(u(t),K(t))$ is elastic-Griffith stable with body load $f(t)$ and boundary load $g(t)$ in the sense of the above Definition \ref{maindefstaticbody}, 
\item $K(t)\subset K(s)$ for all $t<s$, and
\item the total energy,
\[E(t):=\frac12\io |\nabla u(t)|^2 dx -\io f(t)u(t) dx-\int_{\nb} g(t)u(t) \;ds+G_c\huno(K(t))-\mbox{work}(t)\]
is non-increasing in $t$, where work$(t)$ is the work done by applied forces from time zero to time $t$ (as well as by possibly varying Dirichlet boundary values).

\end{enumerate}

\end{defn}
All quasi-static evolutions (without loads) that have been studied mathematically, whether locally stable as in \cite{dmt} and \cite{kmz}, local minimizers as in \cite{le}, or global minimizers as originally proposed in \cite{fm}, satisfy this definition, with $f=g=0$.

We can also define static crack stability for dynamic fracture.  

\begin{defn}[Static stability in dynamic fracture]
For a dynamic fracture solution $(u,K)$ (e.g., as defined in \cite{ladyn}), we say the crack is stable at time $t$ if $(u(t),K(t))$ is locally Griffith stable.

\end{defn}

Notice that we cannot define stability in this way for dynamic fracture when the crack is growing, since where the crack grows, the form of the singularity at crack tips is different from what it is for static and quasi-static fracture.  This is in contrast to a static crack in dynamic fracture, when the form of the stress singularities is the same as it is for static and quasi-static fracture (see \cite{gr2}).  We note that this static stability has not been proven for dynamic fracture solutions.

We now return to the difficulty presented in \cite{sl}.  There, the elastic energy density $W(\nabla u)$ has the form $W^{Eq}(\nabla u) + W^{NEq}(\nabla u)$, and in their model, only $W^{Eq}(\nabla u)$ competes with fracture energy -- that is, cracks advance only if concentrations in $W^{Eq}(\nabla u)$ are sufficiently large.  In terms of energy minimization, this is not so easy to model.  However, with our proposed stability criterion, it is straightforward:  the condition is that every blow-up limit $(\hat{u},\hat{K})$ minimizes
\[E_r(w,\kappa) := \int_{B(0,r)}W^{Eq}(\nabla w) dx +G_c\huno(\kappa \cap B(0,r))
\]
over pairs $(w,\kappa)$ satisfying $\kappa \supset \kh \cap B(0,r)$ and $w\in H_{\hat{u}}^1(B(0,r)\setminus \kappa)$,  for every $r>0$.  So we have

\begin{defn}[Partial-elastic Griffith stability]\label{partial}

$(u,K)$ is $W^{Eq}$-elastic-Griffith stable  if:
\begin{enumerate}
\item $u$ minimizes
\[ w\mapsto \io W(\nabla w) dx\]
over $w\in H_{u}^1(\Omega\setminus K)$, 
and

\item all blow-up limits are $W^{Eq}$-Griffith-stable: for every $\xo\in \Omega$, every blow-up limit $(\hat{u},\hat{K})$ is globally Griffith-stable in the sense of (a modified) Definition \ref{gs}: for every $r>0$, it minimizes
\[E_r(w,\kappa) :=  \int_{B(0,r)}W^{Eq}(\nabla w) dx +G_c\huno(\kappa \cap B(0,r))
\]
over pairs $(w,\kappa)$ satisfying $\kappa \supset \kh\cap B(0,r)$ and $w\in H_{\uh}^1(B(0,r)\setminus \kappa)$.  

\end{enumerate}

In this formulation, we assume that $W^{Eq}$ is quadratic.
Notice also that there is a consistency issue in this definition, not present with the other formulations, or with the phase-field version: if $(u,K)$ satisfies (1), then blow-up limits $(\hat{u},\kh)$ should automatically minimize $E_r(w,\kh)$ over $w\in H_{\uh}^1(B(0,r)\setminus \kh)$.  That is, the blow-up version of the PDE satisfied by $u$ should be the same as the PDE corresponding to this minimization of $E_r$.

\end{defn}

\section{Comparison with stress intensity factors}\label{stress}

A number of mathematics papers have looked carefully at calculating energy release rates using stress intensity factors for regular crack geometries.  In particular, \cite{lt} shows that for sufficiently regular cracks $K$ with crack tip, (antiplane) displacements in elastic equilibrium (with or without loads) must have the form
\[
u(r,\theta)=u_R(r,\theta)+u_S(r,\theta),
\]
where $u_R\in H^2(\Omega\setminus K)$ and $u_S(r,\theta)=\mathcal{K} r^{\frac12}\sin(\theta/2)$ for some $\mathcal{K}\in\R$ (the stress intensity factor), assuming without loss of generality that the crack $K$ is oriented along the $-x$ axis, with crack tip at the origin.  Using \cite{gr2}, the energy release rate is shown to be $\frac\pi4\mathcal{K}^2$, giving crack stability if  $\frac\pi4\mathcal{K}^2< G_c$, and instability if $\frac\pi4\mathcal{K}^2> G_c$.

We now look at what happens with our rescaling in this regular setting: consider $x_0=0$, and we have
\[u_\e(r,\theta):=\e^{-\frac12}u(\e r,\theta)=\e^{-\frac12}u_R(\e r,\theta)+\e^{-\frac12}u_S(\e r,\theta),\] 
where we can assume $u(0)=0$ since that value is subtracted off in the rescaling \eqref{bu}.  The rescaled $u_R$ converge to a constant (since $\nabla u_R$ is bounded - see section 6.4.2 of 
\cite{gr2}), which must be zero since $u_R(0)=0$.  Rescaling $u_S$ gives
\[
(u_S)_\e(r,\theta):=\e^{-\frac12}u_S(\e r,\theta)=\e^{-\frac12}\mathcal{K} (\e r)^{\frac12}\sin(\theta/2)=\mathcal{K} r^{\frac12}\sin(\theta/2)=u_S(r,\theta).
\]
Not surprisingly, we get that the singular part is self-similar under this rescaling.  
We therefore have
\[
u_\e \rightarrow u_S
\]
as $\e\rightarrow 0$, and so  $\uh=u_S$.  

As we now discuss, stability in terms of the energy release rate based on the stress intensity factor $\mathcal{K}$ is consistent with local Griffith stability.

\begin{remark}
Suppose first that we have $(u,K)$ as above, with $\frac\pi4\mathcal{K}^2> G_c$, i.e., the elastic energy release rate exceeds $G_c$.  Performing the rescaling and taking the blow-up limit, we obtain $(\uh,\kh)$ with $\uh(r,\theta)=\mathcal{K} r^{\frac12}\sin(\theta/2)$, as we just described above.  Since $\frac\pi4\mathcal{K}^2> G_c$ by assumption, we get that $(\uh,\kh)$ is not globally Griffith stable, as small increments $\Delta K$ in $K$ would reduce the elastic energy by approximately $\frac\pi4\mathcal{K}^2\huno(\Delta K)> G_c\huno(\Delta K)$, and so such increments would reduce the total energy of $(\uh,\kh)$.  So, instability in the sense that the stress intensity factor is too large implies instability in the sense of local Griffith stability.

On the other hand, suppose we have $(u,K)$ for which $\frac\pi4\mathcal{K}^2< G_c$.  Then we again perform a blow-up, and get $\uh(r,\theta)=\mathcal{K} r^{\frac12}\sin(\theta/2)$.  Since $\frac\pi4\mathcal{K}^2< G_c$, small increments in the crack would increase the total energy, since the release of elastic energy is less than the cost in surface energy.  So, $(\uh,\kh)$ is globally stable, when considering sufficiently small (regular) increments.  But any (regular) increment could be scaled down with the scaling $K_\ep$, making it sufficiently small, without changing $\uh=u_S$ since $u_S$ is self-similar under these rescalings.  So, we get global stability when considering regular increments.  It is a natural conjecture that this minimality extends to all increases in the crack set, but as far as we know, this is open.  

Another question is what happens if $\frac\pi4\mathcal{K}^2=G_c$?  As far as we know, it is open whether $u_S$ in this case would be globally stable, but it is natural to expect it would be.
\end{remark}

\begin{remark}
As noted above, the representation  $u=u_R+u_S$, necessary for the definition of stress intensity factor,  requires some restrictions on $K$, in particular, that $K$ is a connected set near a crack tip.  Our formulation has no such requirement.  For example, $K$ might have a component with crack tip at a point $x_0$, and there could be a family of smaller and smaller cracks approaching $x_0$ from ahead of the crack component with tip $x_0$.  This would prevent the analysis using stress intensity factors and the above representation of $u$ near the crack tip, but our formulation, local Griffith stability, would still be valid.  Indeed, $\uh$ and $\kh$ always exist, and their minimality can always be tested. 

\end{remark}

\begin{remark}
We return now to the question of whether stability should be defined in terms of competition between only elastic energy and surface energy, or whether potential energy due to loads should also be involved.  With sufficient regularity of $K$, the equilibrium displacement has the form $u=u_R+u_S$ with $u_S=\mathcal{K} r^{\frac12}\sin(\theta/2)$, as above.  \cite{gr2} shows that the energy release rate for
\[
\frac12 \int |\nabla u|^2-\int fu
\]
is also $\frac\pi4\mathcal{K}^2$.  But following our definition, we blow-up $u$ at the origin (the crack tip) to get $\uh=u_S$, and we get the same energy release rate of $\frac\pi4\mathcal{K}^2$, even though that energy release rate is based only on $\frac12 \int |\nabla u|^2$.  So as with phase-field fracture, we get the desirable property that stability depends on the stress singularity at a crack tip, and is {\em independent of the reason for the singularity} -- the presence of loads, etc., is irrelevant.
\end{remark}

\begin{remark}
To summarize, taking the blow-up limit is a way of isolating the singular part of $u$, {\em even with no regularity of $K$ and no decomposition of $u$}; $\uh$ is a weak, more general, version of $u_S$, which behaves just like $u_S$ for purposes of stability.

Finally, we emphasize again that this variational formulation works in any dimension, and can accommodate complicated geometries such as those that might be present along crack fronts in three dimensions.

\end{remark}

\section{Minimality and local Griffith stability}\label{mgs}

We now explain why unilateral  global or local minimizers of $E^G$ must be locally elastic-Griffith stable.  Of course, this implies in particular that quasi-static evolutions in the sense of \cite{fm}, and  locally minimizing evolutions, as in \cite{le}, must be locally elastic-Griffith stable at every time.  We also note that if evolutions are not  local minimizers, but are stable in the sense of Griffith at every time, as in \cite{dmt} and \cite{kmz}, then they also satisfy local elastic-Griffith stability at every time.

\subsection{Global Griffith stability and blow-ups}

\subsubsection{Global minimizers}\label{gm}
If $(u,K)$ is a global minimizer of $E^G$ (or unilateral minimizer), then it also satisfies two local properties: i) a weak form of the Euler-Lagrange equation, $\Delta u=0$ in $\Omega\setminus K$ and $\partial_\nu u=0$ on $K$, and ii) local Griffith stability.  Here we outline why ii) holds. First, we restate the scaling -- fix $x_0\in \Omega$, and set 
\begin{equation}
u_\e(x):=\ep^{-\frac12}[u(x_0+\ep x)-u(x_0)]\\
\mbox{ and }K_\e:=\e^{-1}(K-x_0).
\end{equation}
The point of this scaling, introduced in \cite{bo} (and appropriate for any 2-homogeneous elastic energy) is that it preserves the relative contributions of the elastic and surface energies:
\[
\frac12 \int_{B(0,r)}|\nabla u_\e|^2dx = \e^{-1} \frac12  \int_{B(x_0,\e r)}|\nabla u|^2dx
\]
while
\[
\huno(K_\e \cap B(0,r))=\e^{-1} \huno(K\cap B(x_0,\e r)),
\]
with $u_\e$ continuous off of $K_\e$.
For $\Omega\subset \rn$ instead of $\rt$, the only difference is that $\e^{-1} $ is replaced with $\e^{1-N}$ in both of the above equations, and $\huno$ is replaced with $\hn$ in the second.  There is one additional point we should note, which is that more generally, in the definition of $u_\e$, we can subtract off piecewise constant functions (see \cite{bo}) instead of the constant $u(\xo)$.  This is useful to study solutions near points on the crack set besides crack tips.  But for crack tips, it is enough to subtract off $u(\xo)$.

Since the $E^G$ energy of competitors scales the same way, it is immediate that, since by assumption $(u,K)$ has lower energy (elastic plus surface) than competitors that differ from $(u,K)$ only in $B(x_0,\ep r)$, it must be that $(u_\e,K_\e)$ also minimizes this energy (the factor $\e^{-1}$ not affecting minimality) compared with competitors that differ from it only in $B(0,r)$.  To be more explicit, for every $r>0$, given $(w,\kappa)$ such that $w\in H^1_{u_\e}(B(0,r)\setminus \kappa)$ and $\kappa\subset B(0,r)$ with $\kappa\supset K_\e \cap B(0,r)$, we have
\[
\frac12 \int_{B(0,r)}|\nabla u_\e|^2dx +G_c\huno(K_\e \cap B(0,r))\leq \frac12 \int_{B(0,r)}|\nabla w|^2dx +G_c\huno(\kappa).
\]
Notice that this implies 
\be\label{bound}
\frac12 \int_{B(0,r)}|\nabla u_\e|^2dx \leq G_c 2\pi r,
\ee
since $\partial B(0,r)$ can be added to $K_\e \cap B(0,r)$, creating a competitor with zero elastic energy in $B(0,r)$ (technically, $\partial B(0,\rho)$ would be added, with $\rho$ slightly less than $r$).  So, it must be that the amount of elastic energy inside $B(0,r)$ is no more than the perimeter of $B(0,r)$ times $G_c$.

Then, if  $u_\e$ and $K_\e$ converge  to $\hat{u}$ and $\kh$ respectively, as $\e\rightarrow 0$ (or for a sequence $\e_n\rightarrow 0$), it follows from \cite{fl} that for every $r>0$, $(\hat{u},\kh)$ minimizes

\[
(w,\kappa)\mapsto \frac12 \int_{B(0,r)}|\nabla w|^2dx + G_c \huno(\kappa\cap B(0,r))
\]
over $w\in H_{\hat u}^1(B(0,r)\setminus \kappa)$ and $\kappa\supset \kh$.  The idea is that if a pair $(w,\kappa)$ had lower energy than $(\hat{u},\hat{K})$ in $B(0,r)$, then for $\e$ small enough (so $u_\e$ is close enough to $\hat{u}$), $u_\e$ and $K_\e$ could be modified inside $B(0,r)$ to equal $(w,\kappa)$ (at an energy cost that goes to zero with $\e$), thereby reducing the energy of $(u_\e,K_\e)$ in $B(0,r)$, contradicting minimality.  So, these blow-ups satisfy Definition \ref{gs}, and $(u,K)$ is locally elastic-Griffith stable.

\begin{rmk}[On convergence]\label{rmk}
Because of the bound \eqref{bound}, if for each $r>0$, $\huno(K_\e \cap B(0,r))$ is bounded, we are guaranteed that there is some $(\uh,\kh)$ that is the limit of a sequence $(u_{\e_n},K_{\e_n})$ with $\e_n\rightarrow 0$ (this follows from $SBV$ compactness).  The convergence is like weak local $H^1$ convergence, in that $u_\e\rightarrow \uh$ in $L^2(B(0,r))$ and $\nabla u_\e\wc \nabla \uh$ weakly in $L^2(B(0,r))$, for every $r>0$.  If $K_{\e_n}$ are the discontinuity sets of $u_{\e_n}$, then $\kh$ is by definition the discontinuity set of the limit $\uh$.  This convergence is summarized as $(u_{\e_n},K_{\e_n})\rightarrow (\hat{u},\kh)$, and we call $(\uh,\kh)$ a blow-up limit of $(u_{\e},K_{\e})$.
\end{rmk}

\subsubsection{Local minimizers}\label{lm}
It is straightforward to see that blow-ups of local minimizers also satisfy Definition \ref{gs}:   If $(u,K)$ is a ($L^2$) local minimizer, that means that for some $\delta>0$, $(u,K)$ is minimal over competitors $(w,\kappa)$ with $\|u-w\|_{L^2}<\delta$.  Now again consider a blow-up limit $(\hat{u},\hat{K})$.  If $(w,\kappa)$ had lower energy than $(\hat{u},\hat{K})$ in $B(0,r)$, then for $\e$ small enough, the unscaled $w$ would be within $\delta$ of $u$, and $u_\e$ and $K_\e$ could be modified inside $B(0,r)$ to equal $(w,\kappa)$ (at an energy cost that goes to zero with $\e$), thereby reducing the energy of $(u_\e,K_\e)$ in $B(0,r)$. For $\e$ small enough, the $L^2$ norm of this modification would be less than $\delta$.  So, these blow-ups satisfy Definition \ref{gs}.

\subsubsection{Blow-ups of local minimizers with loads}\label{locload}
Whether or not local minimality of $E_{Load}$ is a reasonable model for body loads, the fact is, body loads cannot compete with surface energy locally (on sufficiently small scales), and we can see that here too, minimizers satisfy local elastic-Griffith stability. We now sketch the reason.  We consider $x_0\in\Om$, and for simplicity, we take $u(x_0)=0$.  Then set 
\[u_\e(x):=\ep^{-\frac12}u(x_0+\ep x)\]
 and $K_\e:=\e^{-1}(K-x_0)$.

As noted, this scaling preserves the relative contributions of the elastic and surface energies:
\[
\frac12 \int_{B(0,r)}|\nabla u_\e|^2dx = \e^{-1} \frac12  \int_{B(x_0,\e r)}|\nabla u|^2dx
\]
while
\[
\huno(K_\e \cap B(0,r))=\e^{-1} \huno(K\cap B(x_0,\e r)).
\]

On the other hand, the load term scales differently:
\[
\int_{B(0,r)} f u_\e dx= \e^{-\frac52}\int_{B(x_0,\e r)} f u  dx
\]
so that
\[
\e^{\frac32}\int_{B(0,r)} f u_\e dx = \e^{-1}\int_{B(\xo,\e r)} f u dx.
\]
We then get
\begin{align*}
&E_r(u_\e,K_\e):=\frac12 \int_{B(0,r)}|\nabla u_\e|^2dx -\e^{\frac32}\int_{B(0,r)} f u_\e dx+ \huno(K_\e \cap B(0,r))\\
&= \e^{-1} \frac12  \int_{B(x_0,\e r)}|\nabla u|^2dx- \e^{-1}\int_{B(\xo,\e r)} f u dx+\e^{-1} \huno(K\cap B(x_0,\e r)).
\end{align*}

So, since $u$ was a local minimizer, $u_\e$ minimizes $E_r$ subject to small variations, and so $u_\e$ weakly satisfy 
\[ \Delta u_\e=\e^{\frac32}f.\]
If we suppose  $u_\e$ and $K_\e$ converge  to $\hat{u}$ and $\kh$ as $\e\rightarrow 0$ (at least for some sequence $\e_n\rightarrow 0$), then since $\Delta u_\e\rightarrow 0$, we get
\[ \Delta \hat{u}=0,\]
and we are left with $(\hat u, \kh)$ minimizing just elastic plus surface energy, and so this pair satisfies Definition \ref{gs}.

\section{$\Gamma$-convergence of rescaled phase-field fracture}\label{gcpf}

We now return to the phase-field formulation of fracture, as implemented.  We recall that in practice, the phase-field energies $E_\delta$ are not globally minimized, but perhaps they are  locally minimized.  In this case, $\Gamma$-convergence is not the relevant convergence of phase-field energies to Griffith fracture energies of the type in \cite{fm}.  This raises important questions, such as, what are implementations of phase-field fracture approximating (if anything), and in what sense?  Below we give partial answers.

Suppose $(u_\delta,v_\delta)$ locally minimize $E_\delta[\Om]$, where
\[
E_\delta[\Om](u,v):=\frac12\io (\eta_\delta + v^2)|\nabla u|^2 dx + G_c\left( \frac1{4\delta} \io (1-v)^2 dx + \delta \io |\nabla v|^2 dx  \right),
\]
subject to boundary conditions,
and $u_\delta\rightarrow u_0$.  We want to check whether $(u_0,S_{u_0})$ satisfies elastic-Griffith stability, our Definition \ref{maindefstatic}, wheere $S_{u_0}$ is the discontinuity set of $u_0$.  

For each $u,v\in H^1(\Omega)$ and any $x_0 \in \Om$,  $\e>0$, set
\[u_\e(x):=\ep^{-\frac12}[u(x_0+\ep x)-u(x_0)]\]
 and 
 \[v_\e(x):=v(x_0+\ep x).\]
 Note that $v_\e$ is only a dilation, just as $K_\e$ was a dilation.  Then $|\nabla u_\e|^2(x)=\e|\nabla u|^2(\xo+\e x)$ and $|\nabla v_\e|^2(x)=\e^2|\nabla v|^2(\xo+\e x)$, while scaling from $B(x_0,\e r)$ to $B(0, r)$ produces a factor of $\e^{-2}$, so that
 \[
\frac12 \int_{B(0,r)} (\eta_\delta + v_\e^2)|\nabla u_\e|^2dx = \e^{-1} \frac12  \int_{B(x_0,\e r)}(\eta_\delta + v^2)|\nabla u|^2dx,
\]
\[
\frac1{4\delta} \int_{B(0,r)} (1-v_\e)^2 dx  =\frac{\e^{-2}}{4\delta} \int_{B(x_0,\e r)} (1-v)^2 dx  ,
\]
and 
\[
 \delta \int_{B(0,r)} |\nabla v_\e|^2 dx  = \delta \int_{B(x_0,\e r )} |\nabla v|^2 dx  .
\]
We then get
\[
\frac12 \int_{B(0,r)} (\eta_\delta + v_\e^2)|\nabla u_\e|^2dx+G_c\left[\frac\e{4\delta} \int_{B(0,r)} (1-v_\e)^2 dx + \frac\delta\e \int_{B(0,r)} |\nabla v_\e|^2 dx\right] \]
\[
=\e^{-1}\left(\frac12  \int_{B(x_0,\e r)}(\eta_\delta + v^2)|\nabla u|^2dx+G_c\left[ \frac{1}{4\delta} \int_{B(x_0,\e r)} (1-v)^2 dx + \delta \int_{B(x_0,\e r )} |\nabla v|^2 dx\right]
\right).
\]
Setting $\alpha:=\delta/\e$, we have
\[E_\alpha[B(0,r)](u_\e,v_\e)=\e^{-1}E_\delta[B(x_0,\e r)](u,v).\]
Here we are cheating slightly, since we have $\eta_\delta$ instead of $\eta_\alpha$ in $E_\alpha$, but this makes no difference
 if we choose $\eta_\delta\ll \alpha$.  So we have the usual $\Gamma$-convergence of Ambrosio-Tortorelli \cite{at}, meaning that $E_\alpha[B(0,r)]$ $\Gamma$-converges to the Griffith energy $E^G_r$ for every $r>0$.

Returning to our local minimizers $(u_\delta,v_\delta)$, due to compactness we can suppose that $u_\delta\rightarrow u_0$ in $L^2(\Om)$, and $u_0\in H^1(\Om\setminus K)$ with $K=S_{u_0}$ ($u_0$ is technically in SBV). We fix $x_0\in\Om$ and perform the above $\ep$ rescaling, to get $(u_\delta^\ep, v_\delta^\ep)$, and rescale $u_0$ to get ${u_0}^\e$.  If the blow-up limit of $u_0$ at $x_0$ is $\uh$, i.e., ${u_0}^\e\rightarrow \uh$, then we have
 \[{\displaystyle {\uh= \lim_{\ep\rightarrow 0}{u_0}^\ep}}\mbox{ and }{u_0}^\ep={\displaystyle{ \lim_{\delta\rightarrow 0} u_\delta^\ep}},\] so we get
 \[\uh={\displaystyle{\lim_{\ep\rightarrow 0}\lim_{\delta\rightarrow 0} u_\delta^\ep}}.\]
 We can then choose $\ep(\delta)\gg \delta$ so that ${\displaystyle {\lim_{\delta\rightarrow 0} u_\delta^{\ep(\delta)}=\uh.}}$
 
  If the $(u_\delta,v_\delta)$ are uniform local minimizers near $\xo$, i.e.,   they are all global minimizers within $B(\xo, \rho)$ for some $\rho>0$ (or within $B(\xo,\rho_\delta)$ with $\rho_\delta\rightarrow 0$ more slowly than $\e(\delta)$), then it follows from the $\Gamma$-convergence of $E_\alpha$ that $(\uh,\kh)$ is globally Griffith stable.  And if this uniform local minimality holds for each $\xo\in\Omega$, then $(u_0,S_{u_0})$ satisfies our Definition \ref{lgs}.  Condition (1) in Definition \ref{maindefstatic} is straightforward, and so the pair satisfies Definition \ref{maindefstatic}, and is elastic-Griffith stable.

\section{Conclusions}

We introduce a new, local, variational principle for Griffith fracture which does not have the drawbacks of global minimization, and in particular is compatible with applied loads.  We also emphasize a larger point, that it is not entirely global minimization that leads to problems with loads.  Rather, it is also the requirement that there be a variational formulation using a single energy, so that loads and fracture are forced to directly interact.  Instead, as first put forward with boundary loads in \cite{la}, we should consider two distinct variational problems, one for determining the equilibrium displacement, and the other for defining crack stability based on elastic energy in competition with surface energy.  

Here we showed that with body loads, and in other settings, the competition between elastic and surface energy must be local.   The difficulty is defining a local principle in a way that maintains the mathematical robustness of \cite{fm}.  Local Griffith stability does just this.  In addition, it is consistent with, and can be viewed as a generalization of, stress intensity factors and corresponding energy release rates.  This local stability condition also opens the door to mathematical models for a host of materials which are not purely elastic, but crack growth is governed by elastic energy concentrations competing with fracture surface energy.

\vspace{1cm}

\noindent {\bf Acknowledgments.} This material is based on work supported by the National Science Foundation under Grant No.  DMS-2206114.  The author would also like to thank Oscar Lopez-Pamies for helpful comments on earlier drafts.

\vspace{1cm}

{
}

\end{document}